%Paper: alg-geom/9401002
%From: Joerg.Winkelmann@RUBA.RZ.ruhr-uni-bochum.de
%Date: Sun, 16 Jan 1994 19:48:50 +0100
%Date (revised): Tue, 18 Jan 1994 12:02:39 +0100

%
%  Holomorphic Functions on an Algebraic Group
%  Invariant under Zariski-dense Subgroups
%
%       by Joerg Winkelmann
%
%  WINKELMANN@RUBA.RZ.RUHR-UNI-BOCHUM.DE
%
%
% Lots of private macros
%
\magnification=\magstep1
\catcode`\@=11   %
\def\newfam{\alloc@8\fam\chardef\sixt@@n}
\catcode`\@=12
\def\today{January 15, 1994}
\font\petcap=cmcsc10
\font\sevenit=cmti7 scaled \magstep0
\font\sevensl=cmti7 scaled \magstep0
\font\bigit=cmti10 scaled \magstep1
\font\bigbf=cmbx10 scaled \magstep1
\font\bigrm=cmr10 scaled \magstep1
\font\bigsevenbf=cmbx7 scaled \magstep1
\font\bigsevenrm=cmr7 scaled \magstep1
\font\bigteni=cmmi10 scaled \magstep 1
\font\bigex=cmex10 scaled \magstep1
\font\bigsy=cmsy10 scaled \magstep1
\font\bigseveni=cmmi7 scaled \magstep1
\font\bigsevensy=cmsy7 scaled \magstep1
\skewchar\bigteni='177
\skewchar\bigsy='60
\catcode`\@=11
\newskip\standardbaselineskip
\standardbaselineskip=14pt
\def\tenpoint{%
 \textfont0=\tenrm
 \textfont1=\teni
 \textfont2=\tensy
 \textfont\bffam=\tenbf
 \textfont\ttfam=\tentt
 \textfont\itfam=\tenit
 \textfont\slfam=\tensl
\scriptfont0=\sevenrm \scriptscriptfont0=\fiverm
 \scriptfont1=\seveni \scriptscriptfont1=\fivei
 \scriptfont2=\sevensy \scriptscriptfont2=\fivesy
 \scriptfont3=\tenex \scriptscriptfont3=\tenex
 \scriptfont\bffam=\sevenbf \scriptscriptfont\bffam=\fivebf
\scriptfont4=\sevenit
\scriptfont5=\sevensl
  \def\it{\fam\itfam\tenit}%
  \def\rm{\fam\z@\tenrm}%
  \def\bf{\fam\bffam\tenbf}%
  \def\sl{\fam\slfam\tensl}%
  \def\tt{\fam\ttfam\tentt}%
  \def\oldstyle{\fam\@ne\teni}%
  \def\big##1{{\hbox{$\left##1\vbox to 8.5pt{}\right.\n@space$}}}
  \setbox\strutbox=\hbox{\vrule height8.5pt depth3.5pt width0pt}%
  \abovedisplayskip=13pt plus 3pt minus 10pt
  \belowdisplayskip=13pt plus 3pt minus 10pt
  \normalbaselineskip=\standardbaselineskip
  \normalbaselines\rm}
\def\ninepoint{\tenpoint}
\def\bigstyle{
\textfont0=\bigrm
\scriptfont0=\bigsevenrm
\textfont1=\bigteni
\scriptfont1=\bigseveni
\textfont2=\bigsy
\scriptfont2=\bigsevensy
\textfont3=\bigex
\textfont4=\bigit
\scriptfont4=\tenit
\textfont5=\bigit %
\textfont6=\bigbf
\scriptfont6=\bigsevenbf
\textfont\ttfam=\nullfont
  \def\it{\fam\itfam\bigit}%
  \def\rm{\fam\z@\bigrm}%
  \def\bf{\fam\bffam\bigbf}%
  \def\sl{\fam\slfam\bigit}%
  \def\tt{\errmessage{\string\tt\ in \string\bigstyle\ not yet implemented}}%
  \def\oldstyle{\fam\@ne\bigteni}%
  \def\big##1{{\hbox{$\left##1\vbox to 10.2pt{}\right.\n@space$}}}
  \setbox\strutbox=\hbox{\vrule height10.2pt depth4.2pt width0pt}%
  \normalbaselineskip=\standardbaselineskip
  \advance\normalbaselineskip by 3pt
  \normalbaselines\rm}
\catcode`\@=12
\tenpoint
\def\Bbb#1{{\bf #1}}

\def\eg{e.g.\spacefactor=1000\ }%
\def\ie{i.e.\spacefactor=1000\ }%
\let\oldcal\cal
\def\cal#1{{\oldcal #1}}%
\let\Cedille\c
\def\abs#1{{\vert #1\vert}}

\def\sdr{\mathbin{\vrule height 4.5pt depth 0.1pt\kern-1.8pt\times}}
\def\={\hbox{--}\nobreak\hskip0pt\relax}
\def\({\left(}
\def\){\right)}%
\pageno=1
\relpenalty=5000
\binoppenalty=5000
\hfuzz=2truept
\mathsurround=1pt
\nulldelimiterspace=0pt
\overfullrule=0pt
\normallineskip=1pt
\normallineskiplimit= 1pt
\normalbaselines
\newdimen\oldparindent

\def\restoreparindent{\let\parindent\oldparindent}%
\def\.{\null.}%
\newtoks\autor
\autor={J\"org Winkelmann}%
\newtoks\ort
\ort={Bochum}%
\newcount\arbvar
\newtoks\normalheadline
\normalheadline={\hfil}
\headline={\ifnum\pageno=1\hfill\else\the\normalheadline\fi}%
\catcode`\@=11
\def\makeheadline{\vbox to\z@{\vskip-32.5\p@
  \line{\vbox to8.5\p@{}\the\headline}\vss}\nointerlineskip}%
\catcode`\@=12
\newtoks\normalfootline
\footline={\hfill\tenrm\folio\hfill}

\newskip\vardimen
\newskip\varlastdimen
\def\varskip{\varvarskip{0pt}}%
\def\varvarskip#1#2\penalty#3 {\varlastdimen=\lastskip
\removelastskip\penalty#3\relax
\vardimen=#2\relax
\advance\varlastdimen by#1\relax
\ifdim\vardimen<\varlastdimen\relax\vardimen=\varlastdimen\fi
\vskip\vardimen\relax}%
\catcode`\@=11
\def\footnote#1{\let\@sf\empty %
  \ifhmode\edef\@sf{\spacefactor\the\spacefactor}\/\fi
 $^{\rm #1}$\@sf\vfootnote{#1}}
\def\vfootnote#1{\insert\footins\bgroup
  \interlinepenalty\interfootnotelinepenalty
  \ninepoint
  \splittopskip\ht\strutbox %
  \splitmaxdepth\dp\strutbox \floatingpenalty\@MM %
  \leftskip\z@skip \rightskip\z@skip \spaceskip\z@skip \xspaceskip\z@skip
  \textindent{$^{\rm #1}$}\footstrut\futurelet\next\fo@t}%
\catcode`\@=12
\newtoks\postdisplaytoks
\everydisplay=\expandafter{\the\everydisplay\global\postdisplaytoks={}%
\aftergroup\thepostdisplaytoks}%
\def\thepostdisplaytoks{\the\postdisplaytoks}%
\def\afterdisplay#1{\global\postdisplaytoks=\expandafter{\the\postdisplaytoks
#1}}%
\newif\ifdisplay
\global\displayfalse
\everydisplay=\expandafter{\the\everydisplay\displaytrue}%
\newcount\qedpenalty
\def\qed{\qedsign\ifdisplay\postdisplaypenalty=\qedpenalty
\afterdisplay\endproofstyle
\let\next\relax\else\par\penalty\qedpenalty\let\next\endproofstyle\fi\next}%
\def\qedsign{\ifmmode\ \qedkasten{4pt}\else\nobreak\ $\qedkasten{4pt}$\fi}%
\def\qedkasten#1{\mathop{\kern0.5pt\vbox{\hrule\hbox{\vrule%
   \hskip#1\vrule height #1 width0pt%
   \vrule}\hrule}\kern-0.5pt}}%
\qedpenalty=-100
\def\endproofstyle{\tenpoint\medskip}%
\def\s#1{\Bbb #1}%
\normalheadline=%
{\ninepoint\sl\hfil\ifodd\pageno\the\shorttitle
\else\the\autor\fi\hfil}%
\newbox\titbox
\def\fintitre{\global\let\titre\relax}%
\def\titel#1{\bigskip\bigskip
\def\titre##1\\{
\setbox\titbox=%
\centerline{\bigstyle\bf ##1}%
\ifdim\ht\titbox>0pt\relax \box\titbox\medskip\fi
\titre}%
\titre#1\fintitre\\
\medskip
\centerline{\sl by}
\medskip
\centerline{\rm\the\autor}
\bigskip\bigskip
}%
\newtoks\shorttitle
\shorttitle={}%
\newtoks\autoraddress
\autoraddress={%
J\"org Winkelmann\par
Mathematisches Institut\par
NA 4/69 \par
Ruhr-Universit\"at Bochum \par
D-44780 Bochum \par
Germany \par
\medskip
e-mail (BITNET):\par
WINKELMANN\par
@RUBA.RZ.RUHR-UNI-BOCHUM.DE\par
}%
\def\datum{\par\medskip\bgroup\leftskip=200pt \parindent=0pt \the\ort,
\today\par\bigskip\egroup}
\def\endarticle{\penalty100\relax\bigskip\datum\penalty-300
\smallskip %
\theaddress\end}%
\def\theaddress{%
\vbox{\leftskip=200pt
\parindent=0pt %
\baselineskip=8.5pt
\interlinepenalty=4000
\rm\the\autoraddress}}
\newtoks\thankstoks
\def\thanks#1{\thankstoks={#1}}
\def\thethanks{\par\the\thankstoks\par}
\newcount\sectno
\sectno=0
\def\section#1{\varskip\bigskipamount\penalty-250
\advance\sectno by 1
\leftline{\ifnum\sectno>0\relax\the\sectno. \else\fi
\petcap\ignorespaces #1}%
\nobreak
\medskip}%
\def\subsection#1{\varskip\medskipamount\penalty-120
\centerline{\petcap #1}\nobreak\medskip}%
\def\subsubsection#1{\varskip\medskipamount\penalty-120
\bgroup\it\ignorespaces#1. \egroup\ignorespaces}
\catcode`\@=11
\bgroup\obeylines\gdef\proclaim@@#1#2#3#4#5
{\egroup\varskip\bigskipamount\penalty30 %
\noindent{#1#2%
\def\parameter{#5}%
\ \ignorespaces\ifx\parameter\empty #3\else #5\fi\unskip. }%
\bgroup#4\ignorespaces}\egroup
\def\endproclaim{\par\egroup\penalty-50\medskip}%
\def\proclaim@{\bgroup\obeylines\proclaim@@}%
\def\proclaim#1{\proclaim@{\bf}{#1}{}{\it}}
\def\Corollary{\proclaim@{\bf}{Corollary}{}{\it}}
\def\Example{\proclaim@{\bf}{Example}{}{\rm}}
\def\Definition{\proclaim@{\bf}{Definition}{}{\rm}}
\def\Remark{\proclaim@{\bf}{Remark}{}{\rm}}
\def\Proof{\proclaim@{\it}{Proof}{}{\ninepoint\rm}}
\newcount\propno
\propno=1
\def\Proposition{\proclaim@{\bf}{Proposition}
{\the\propno\global\advance\propno by 1 }{\it}}
\newcount\theono
\theono=1
\def\Theorem{\proclaim@{\bf}{Theorem}
{\the\theono\global\advance\theono by 1 }{\it}}
\newcount\lemno
\lemno=1
\def\Lemma{\proclaim@{\bf}{Lemma}{\the\lemno\global\advance\lemno by 1 }{\it}}
\catcode`\@=12

\def\endproofstyle{\egroup}%
\newif\ifciteerror
\newcount\citeno
\def\checkciteerror{\ifciteerror\errmessage{Citation error}\fi}
\def\defcitations#1{\dodefcites#1,,\enddefcite}
\def\dodefcites#1,{\def\next{#1}\ifx\next\empty
\let\next\finishdefcite
\else\let\next\doodefcite\fi
\next{#1}}
\def\doodefcite#1{\advance\citeno by 1
\expandafter\xdef\csname cite:#1\endcsname
{\noexpand\global\noexpand\citeerrorfalse
\the\citeno}\dodefcites} \def\finishdefcite#1\enddefcite{}
\def\makecite#1#2{\csname cite:#1\endcsname#2}
\def\docite#1,#2,#3!{\def\next{#2}\ifx\next\empty[\makecite{#1}{}]\else
[\makecite{#1}{,#2}]\fi}
\def\cite#1{{\rm\citeerrortrue\docite#1,,!\checkciteerror}}
\newdimen\refindent
\def\References{\medskip\sectno=-100 \section{References}
\let\c\Cedille
\def\(##1\){}
\frenchspacing
\setbox0=\hbox{99.\enspace}\refindent=\wd0\relax
\def\item##1##2:{\goodbreak
\hangindent\refindent\hangafter=1\noindent\hbox
to\refindent{\hss\refcite{##1}\enspace}\petcap
\ignorespaces##2\rm: }}
\def\refcite#1{\citeerrortrue\dorefcite#1[]\endcite\checkciteerror.}
\def\dorefcite#1[#2]#3\endcite{\def\next{#2}%
\ifx\next\empty\makecite{#1}{}\else\makecite{#2}{}\fi}%
\def\AMS#1{\relax}
%%%%%
%%%%%%
\catcode`\@=11
\def\Examples{\proclaim@{\bf}{Examples}{}{\rm}}
\catcode`\@=12
\def\implies{\Rightarrow}
\def\Lie{\mathop{{\cal L}{\it ie}}}
\def\Gr{G_{\s R}}
\defcitations{BO,BKO,HO,HM,M1,M2,Mi,T}
%%%%%%%%%%%%%%%%%%%%%%%%%%%%%%%
% Here starts the mathematics %
%%%%%%%%%%%%%%%%%%%%%%%%%%%%%%%
\titel{Holomorphic Functions on an Algebraic Group \\
Invariant under Zariski-dense Subgroups}
\section{Introduction}
Let $G$ be a reductive complex linear-algebraic group, $\Gamma$ a subgroup,
$\cal O(G)^\Gamma$ the algebra of $\Gamma$-invariant holomorphic functions
on $G$. It is known \cite{BO} that $\cal O(G)^\Gamma=\s C$ if $\Gamma$ is dense
in $G$ with respect to the algebraic Zariski-topology.

We are interested in similar results for non-reductive groups.
If $G$ is a complex linear-algebraic group with $G/G'$ non-reductive
(where $G'$ denotes the commutator group),
then there exists a surjective group morphism $\tau:G\to(\s C,+)$ and
$\Gamma=\tau^{-1}(\s Z)$ is a Zariski-dense subgroup of $G$ with
$\cal O(G)^\Gamma\simeq\cal O(\s C^*)\ne\s C$.
Hence we are led to the question whether the following two properties
are equivalent:

Let $G$ be a connected complex linear-algebraic group.
\item{(i)} $G/G'$ is reductive.
\item{(ii)} $\cal O(G)^\Gamma=\s C$ for every Zariski-dense subgroup
$\Gamma$.

The above argument gave us $(ii)\implies (i)$ and the result of Barth and
Otte \cite{BO}
implies the equivalence of $(i)$ and $(ii)$ for $G$ reductive.

We will prove that $(i)$ and $(ii)$ are likewise equivalent in the
following two cases:

\item{a)} $G$ is solvable.
\item{b)} The adjoint representation of $S$ on $\Lie(U)$ has no zero
weight, where $S$ denotes a maximal connected semisimple subgroup of $G$
and $U$ the unipotent radical of $G$.

Case $b)$ is equivalent to each of the following two conditions
\item{b')}
$G/G'$ is reductive
and the semisimple elements are dense in $G'$.
\item{b'')}
$G/G'$ is reductive and $N_{G'}(T)/T$ is finite, where $T$ is a maximal
torus in $G'$ and $N_{G'}(T)$ denotes the normalizer of $T$ in $G'$.

For instance, if we take $G$ to be a semi-direct product $SL_2(\s C)\sdr_\rho
(\s C^n,+)$ with $\rho:SL_2(\s C)\to GL_n(\s C)$ irreducible, then
$G$ fulfills the condition of case $b)$ if and only if $n$ is an even
number.

The proof for case $a)$ is based on the usual solvable group methods and
the structure theorem on holomorphically separable solvmanifolds
by Huckleberry and E. Oeljeklaus.

The proof for case $b)$ relies on the discussion of semisimple elements
of infinite order in such a $\Gamma$.
For this reason we conclude the paper with an example of Margulis which
implies that $G=SL_2(\s C)\sdr_\rho (\s C^3,+)$ ($\rho$ irreducible)
admits a Zariski-dense discrete subgroup $\Gamma$
such that no element of $\Gamma$ is semisimple.
Thus condition $b)$ is really needed in order to find semisimple elements
in Zariski-dense subgroups.
In consequence, our method does not work for the example of Margulis.
However, this only means that we can not prove
$\cal O(G)^\Gamma=\s C$ for Margulis' example.
We have no knowledge whether there actually exist non-constant holomorphic
functions in this case.

Finally we discuss invariant meromorphic and plurisubharmonic functions
on certain groups.

\section{Solvable groups}
Here we will discuss solvable groups. First we will develope some
auxiliary lemmata.
\Lemma
Let $G$ be a connected complex linear-algebraic group such that
$G/G'$ is reductive.

Then $[G,G']=G'$.
\endproclaim
\Proof
By taking the appropriate quotient, we may assume $[G,G']=\{e\}$.
We have to show that this implies $G'=\{e\}$.
Now $[G,G']=\{e\}$ means that $G'$ is central, hence $Ad(G)$ factors
through $G/G'$. But $G/G'$ is reductive and acts trivially
(by conjugation) on both $G/G'$ and $G'$. Due to complete
reducibility of representations of reductive groups
it follows that $Ad(G)$ is
trivial, \ie $G$ is abelian, \ie $G'=\{e\}$.
\qed
\Lemma
Let $G$ be a connected complex linear-algebraic group, $H\subset G'$
a connected complex Lie subgroup which is normal in $G'$.

Then $H$ is algebraic.
\endproclaim
\Proof
Let $U$ denote the unipotent radical of $G'$. Then $G'/U$ is semisimple.
Now $A=(H\cap U)^0$ is algebraic, because every connected complex Lie subgroup
of a unipotent group is algebraic. Normality of $H$ implies that
$H/(H\cap U)$ is semisimple. Hence $H/A$ is semisimple, too.
It follows that $H/A$ is an algebraic subgroup of $G'/A$.
Thus $H$ has to be algebraic.
\qed

\Lemma
Let $G$ be a connected topological group, $H$ a normal subgroup,
such that $H\cap G'$ is totally disconnected.

Then $H$ is central.
\endproclaim
\Proof
For each $h\in H$ the set $S_h=\{ghg^{-1}h^{-1}:g\in G\}$ is both
totally disconnected and connected and therefore reduces to $\{e\}$.
\qed

\Lemma
Let $G$ be a connected complex linear-algebraic group,
$A\subset G'$ a complex Lie subgroup which is normal in $G$ and
Zariski-dense in $G'$.
Assume moreover that $[G,G']=G'$.

Then $A=G'$.
\endproclaim
\Proof
The connected component $A^0$ of $A$ is algebraic (Lemma 2).
Thus $G/A^0$ is again algebraic. Moreover $G'/A^0$ is the commutator group
of $G/A^0$.
Therefore, by replacing $G$ with $G/A^0$ we may assume that $A$ is
totally disconnected.
But totally disconnected normal subgroups of connected Lie groups are central
(Lemma 3).
Since $A$ is Zariski-dense in $G'$, and $[G,G']=G'$, this may occur only
for $A=G'=\{e\}$.
Thus $A=G'$.
\qed

\Theorem
Let $G$ be a connected solvable complex linear-algebraic group,
$\Gamma$ a subgroup which is dense in the algebraic
Zariski-topology.
Assume that $G/G'$ is reductive.

Then $\cal O(G)^\Gamma=\s C$.
\endproclaim
\Proof
Let $G/\Gamma \to G/H$ denote the holomorphic reduction,
\ie $$H=\{g\in G: f(g)=f(e)\forall f\in\cal O(G)^\Gamma\}.$$
Now $\Gamma\subset H$ normalizes $H^0$.
Since $\Gamma$ is Zariski-dense in $G$ and the normalizer of a
connected Lie subgroup is necessarily algebraic, it follows that
$H^0$ is normal in $G$.
Let $A=H^0\cap G'$. This is again a closed normal subgroup in $G$.
By a result of Huckleberry and E.\ Oeljeklaus \cite{HO}
$H/H^0$ is almost
nilpotent (\ie admits a subgroup of finite index which is nilpotent).
Let $\Gamma_0$ be a subgroup of finite index in $\Gamma$ with
$\Gamma_0/(\Gamma_0\cap H^0)$ nilpotent.
By definition this means there exists a number $k$ such that
$C^k\Gamma_0\subset H^0$ where $C^k$ denotes the central series.
Now $[G,G']=G'$ implies $C^kG=G'$ for all $k\ge 1$.
Therefore $C^k\Gamma_0$ is Zariski-dense in $G'$.
It follows that $A=H^0\cap G'$ is a closed normal Lie subgroup of $G$
which is Zariski-dense in $G'$.
By the preceding lemma it follows that $A=G'$, \ie $G'\subset H$.
Now $G/G'$ is assumed to be reductive.
Thus the statement of the theorem now follows from the result
for reductive groups (\cite{BO}).
\qed
\section{Groups with many semisimple elements}
Here we will prove the following theorem.
\Theorem
Let $G$ be a connected complex linear-algebraic group.
Assume that $G/G'$ is reductive and that furthermore one
(hence all) of the following equivalent conditions is fulfilled.

\item{(1)}
$G'$ contains a dense open subset $\Omega$ such that each element in
$\Omega$ is semisimple.
\item{(2)}
For any maximal torus $T$ in $G'$ the quotient $N_{G'}(T)/T$
is finite.
\item{(3)}
Let $S$ denote a maximal connected semisimple subgroup of $G$ and $U$ the
unipotent radical of $G$. Let $\rho:S\to GL(\Lie U)$ denote the representation
obtained by restriction from the adjoint representation $Ad:G\to GL(\Lie G)$.
The condition is that all weights of $\rho$ are non-zero.

Under these assumptions $\cal O(G)^\Gamma=\s C$ for any
Zariski-dense subgroup $\Gamma\subset G$.
\endproclaim
\Examples
\item{(a)} Let $G$ be a reductive group. Then $G/G'$ is reductive and
$G'$ semisimple, hence $N_{G'}(T)/T$ finite for any maximal torus
$T\subset G'$. Therefore this theorem is a generalization of the
result of Barth and Otte \cite{BO} on redutive groups.
\item{(b)}
Let $G$ be a parabolic subgroup of a semisimple group $S$.
$G/G'$ is obviously reductive.
Furthermore a maximal torus $T$ in $G$ is already a maximal torus in $S$.
Hence $N_S(T)/T$ is finite. Consequently $N_G(T)/T$ is finite and
$G$ fulfills the assumptions of the theorem.
\item{(c)}
Let $G$ be a semi-direct product of $SL_2(\s C)$ with a unipotent group
$U\simeq\s C^n$ induced by an irreducible representation
$\xi:SL_2(\s C)\to GL(U)$.
Then $G$ fulfills the assumptions of
the theorem if and only if $n$ is even.
\endproclaim

Now we will demonstrate that $(1)$, $(2)$ and $(3)$ are indeed
equivalent.
The equivalence of $(2)$ and $(3)$ is rather obvious from standard results
on algebraic groups.
For the equivalence of $(1)$ and $(2)$ we need some elementary facts
on semisimple elements in a connected algebraic group $G$.
Let $G_s$ denote the set of all semisimple elements in $G$ and $T$ be a
maximal torus in $G$. Now $g\in G_s$ iff $g$ is conjugate to an element in
$T$. It follows that $G_s$ is the image of the map $\zeta:G\times T\to G_s$
given by $\zeta(g,t)=gtg^{-1}$. In particular $G_s$ is a constructible
set. Now a torus contains only countably many algebraic subgroups, hence
a generic element $h\in T$ generates a Zariski-dense subgroup of $T$.
It follows that for a generic element $h\in T$ the assumption $g\in G$
with $ghg^{-1}\in T$ implies $gTg^{-1}=T$.
{}From this it follows that a generic fiber of $\zeta$ has the dimension
$\dim N_G(T)$. Therefore the dimension of $G_s=Image(\zeta)$
equals $\dim G-\dim N_G(T)$. Thus we obtained the following lemma,
which implies the equivalence of $(1)$ and $(2)$.

\Lemma
Let $G$ be a connected linear-algebraic group, $T$ a maximal torus and
$G_s$ the set of semisimple elements in $T$.

Then $G_s$ is dense in $G$ if and only if $\dim N_G(T)=\dim T$.
\endproclaim

Next we state some easy consequences of the assumptions of Theorem 2.
\Lemma
Let $G$ be an algebraic group fulfilling the assumptions of Theorem 2
and $\tau:G\to H$ a surjective morphism of algebraic groups.

Then $H$ likewise fulfills the assumptions of Theorem 2.
\endproclaim
\Proof
Surjectivity of $\tau$ gives a surjective morphism of algebraic
groups from $G/G'$ onto $H/H'$.
Therefore $H/H'$ is reductive.
The surjectivity of $\tau$ furthermore implies $\tau(G')=H'$.
Since morphisms of algebraic groups map semisimple elements to
semisimple elements, it follows that $H$ fulfills condition $(1)$.
\qed
\Lemma
Let $G$ be an algebraic group fulfilling the assumptions of Theorem 2.
Then the center $Z$ of $G$ must be reductive.
\endproclaim
\Proof
Condition $(2)$ implies that $(Z\cap G')^0$ is contained in a maximal
torus of $G'$. Since $G/G'$ is reductive, this implies that $Z$ is reductive.
\qed
The following lemma illuminates why semisimple elements are important
for our purposes.
\Lemma
Let $G$ be a complex linear-algebraic group, $g\in G$ an element of infinite
order, $\Gamma$ the subgroup generated by $g$ and $H$ the Zariski-closure
of $\Gamma$.

Then $Z=H/\Gamma$ is a Cousin group (hence in particular $\cal O(Z)=\s C$)
if $g$ is semisimple; but $Z$ is biholomorphic to some $(\s C^*)^n$ (hence
holomorphically separable) if $g$ is not semisimple.
\endproclaim
\Proof
Note that $\bar\Gamma=H$ implies $H=H^0\Gamma$.
Hence $H/\Gamma=H^0/(H^0\cap\Gamma)$ is connected.
If $g$ is semisimple, the Zariski-closure of $\Gamma $ is reductive
and the statement follows from \cite{BO}.
If $g$ is not semisimple then $H\simeq(\s C^*)^{n-1}\times\s C$ for some
$n\ge 1$ and $g$ is not contained in the maximal torus of $H$.
This implies $H/\Gamma\simeq(\s C^*)^n$.
\qed

%%%%%%%%%%%%%%%%%%
\Lemma
Let $G$ be a connected real Lie group, $\Gamma$ a subgroup such that
each element $\gamma\in\Gamma$ is of finite order.

Then $\Gamma$ is almost abelian and relatively compact in $G$.
\endproclaim
\Proof
If $G$ is abelian, then $G\simeq\s R^k\times(S^1)^n$. In this case
$\Gamma\subset(S^1)^n$ and the statement is
immediate.

Now let us assume that $G$ may be embedded into a complex linear-algebraic
group $\tilde G$. Let $H$ denote the (complex-algebraic) Zariski-closure
of $\Gamma$ in $\tilde G$.
By the theorem of Tits \cite{T} $\Gamma$ is almost solvable,
hence $H^0$ is solvable.
Now the commutator group of $H^0$ is unipotent and therefore contains no
non-trivial element of finite order. Hence $\Gamma\cap H^0$ is abelian,
which completes the proof for this case,
since we discussed already the abelian case.

Finally let us
discuss the general case. By the above considerations $Ad(\Gamma_0)$
is contained in an abelian connected compact subgroup $K$ of $Ad(G)$
for some subgroup $\Gamma_0$ of finite index in $\Gamma$.
Now $N=(Ad)^{-1}(K)$ is a central extension $1\to Z\to N\to K\to 1$.
(where $Z$ is the center of
$G$). But complete reducibility of the
representations of compact groups implies that this sequence splits on the
Lie algebra level. Hence $N$ is abelian and we can complete the proof as
before.
\qed
\Lemma
Let $G$ be a complex linear-algebraic group, $\Gamma$ a Zariski-dense
subgroup.

Then $\Gamma$ contains a finitely generated subgroup $\Gamma_0$ such that
the Zariski-closure of $\Gamma_0$ contains $G'$.
\endproclaim
\Proof
Consider all finitely generated subgroups of $\Gamma$ and their
Zariski-closure in $G$. There is one such group $\Gamma_0$ for which the
dimension of the Zariski-closure $A$ is maximal.
Clearly $A$ must contain the connected component of the Zariski-closure
for any finitely generated subgroup of $\Gamma$.
This implies that $A^0$ is normal in $G$.
Furthermore maximality implies that the group $\Gamma/A^0$ contains no
element of infinite order. Hence $\Gamma/A^0$ is almost abelian, which
implies $G'\subset A$ ($\Gamma/A^0$ is Zariski-dense in $G/A^0$).
\qed
Caveat: There is no hope for $A=G$, \eg take $G=\s C^*$ and let $\Gamma$
denote the subgroup which consists of all roots of unity.

Theorem 2 follows by induction on $dim(G)$ using the following lemma.
\Lemma
Let $G$ be a positive-dimensional
complex linear-algebraic group fulfilling the assumptions of Theorem 2,
$\Gamma$ a Zariski-dense subgroup.

Then there exists a positive-dimensional normal algebraic subgroup $A$
with $\cal O(G)^\Gamma\subset\cal O(G)^A$.
\endproclaim
\Proof
If $G$ is abelian, the assumptions imply that $G$ is reductive and
$\cal O(G)^\Gamma=\s C$.

Otherwise let $H=\{g:f(g)=f(e)\forall f\in\cal O(G)^\Gamma\}$.
Now $\Gamma\subset H$, hence $\Gamma$ normalizes $H^0$.
The normalizer of a connected Lie subgroup is algebraic, thus
$H^0$ is a normal subgroup of $G$.
It follows that $(H\cap G')^0$ is a normal algebraic
subgroup (Lemma 2). This completes the proof unless $H\cap G'$ is
discrete.

This leaves the case where $(H\cap G')$ is discrete.
Then $H^0$ is contained in the center $Z$ of $G$ (Lemma 3).
Let $A$ denote the Zariski-closure of $H^0$.
The center is reductive (Lemma 7). It follows that for each $Z$-orbit every
$H^0$-invariant functions is already $A$-invariant.

Therefore we can restrict to the case where $H$ is discrete.
Now $\Gamma$ is discrete and contains a subgroup $\Gamma_0$ which is
finitely generated and whose Zariski-closure contains $G'$.
By a theorem of Selberg $\Gamma_0$ contains a subgroup of finite index
$\Gamma_1$ which is torsion-free.
Now let $\Gamma_2=\Gamma_1\cap G'$.
Then being Zariski-dense, $\Gamma_2$ must contain a semisimple element
of infinite order. Using Lemma 8, this yields
 a contradiction to the assumption that
$H$ is discrete.
\qed

\section{An example}
At a first glance, it seems to be obvious that a Zariski-dense subgroup
should contain enough elements of infinite order to generate a subgroup
which is still Zariski-dense.
However, one has to careful.
\Lemma
Let
$G=\s C^*\sdr\s C$ with group law $(\lambda,z)\cdot(\mu,w)=(\lambda\mu,
z+\lambda w)$ and $\Gamma$ the subgroup generated by the elements
$a_n=(e^{2\pi i/n},0)$ ($n\in\s N)$ and $a_0=(1,1)$.
Then $\gamma\in G'=\{(1,x):x\in \s C\}$ for any element $\gamma\in\Gamma$
of infinite order, although $\Gamma$ is Zariski-dense in $G$.
\endproclaim
\Proof
It is clear that $\Gamma$ is Zariski-dense in $G$. The other assertion
follows from the fact that any element in $G$ is either unipotent or
semisimple. Hence every element $g\in G\setminus G'$ is conjugate to
an element in $\s C^*\times\{0\}$.
\qed

\section{Margulis' example}
We will use an example of Margulis to demonstrate the following.
\Proposition
There exists a discrete Zariski-dense subgroup $\Gamma$ in
$G=SL_2(\s C)\sdr_\rho(\s C^3,+)$ with $\rho$ irreducible
such that $\Gamma$ contains no semisimple element.
\endproclaim
Thus the condition $G/G'$ reductive is not sufficient to guarantee the
existence of semisimple elements in Zariski-dense subgroups.

Margulis \cite{M1,M2}
constructed his example in order to prove that there exist
free non-commutative groups acting on $\s R^n$ properly discontinuous
and by affin-linear transformations, thereby contradicting a conjecture
of Milnor \cite{Mi}.

We will now start with the description of Margulis' example.
Let $B$ denote the bilinear form on $\s R^3$ given by
$B(x_1,x_2,x_3)=x_1^2+x_2^2-x_3^2$, $W=\{x\in \s R^3:B(x,x)=0\}$ the zero
cone and $W^+=\{x\in W:x_3>0\}$ the positive part.
Let $S=\{x\in W^+:\abs{x}=1\}$.
Let $H$ be the connected component of the isometry group $O(2,1)$ of $B$.
(As a Lie group $H$ is isomorphic to $PSL_2(\s R)$.)
Let $\Gr=H\sdr(\s R^3,+)$ the group of affine-linear transformations
on $\s R^3$ whose linear part is in $H$.

The following is easy to verify:
Let $x^+,x^-$ to different vectors in $S$.
Then there exists a unique vector $x^0$ such that $B(x^0,x^0)=1$,
$B(x^+,x^0)=0=B(x^-,x^0)$ and $x^0,x^-,x^+$ form a positively
oriented basis of the vector space $\s R^3$.
Furthermore for any $\lambda\in ]0,1[$ there is an element $g\in H$
(depending on $x^+,x^-\in S$ and $\lambda$)
defined as follows:
$$ g:ax^0+bx^-+cx^+ \mapsto ax^0+{b\over\lambda}x^-+\lambda c x^+$$
Conversely any non-trivial
diagonalizable element $g\in H$ is given in such a way and
$x^+$, $x^-$ and $\lambda$ are uniquely determined by $g$.

The result of Margulis is the following:

Let $x^+,x^-,\tilde x^+,\tilde x^-$ four different points in $S$,
$\lambda,\tilde\lambda\in]0,1[$, $v,\tilde v\in\s R^3$ such that
$v,x^-,x^+$ resp. $\tilde v,\tilde x^-,\tilde x^+$ forms a positively
oriented basis of $\s R^3$.
Let $h,\tilde h\in H$ be the elements corresponding to $x^-,x^+,\lambda$
resp. $\tilde x^-,\tilde x^+,\tilde\lambda$ and $g,\tilde g\in
\Gr=H\sdr\s R^3$ given by $g=(h,v)$, $\tilde g=(\tilde h,\tilde v)$.

Then there exists a number $N=N(g,\tilde g)$
such that the elements $g^N$, $\tilde g^N$ generate a (non-commutative) free
discrete subgroup $\Gamma\subset\Gr$ such that the action on $\s R^3$ is
properly discontinuous and free.

Now an element $g\in\Gr$ is conjugate to an element in $H$ if and only if
$g(w)=w$ for some $w\in\s R^3$. Hence no element in $\Gamma$ is conjugate
to an element in $H$. In particular no element in semisimple.
Furthermore it is clear that $\Gamma$ is Zariski-dense in the complexification
$G=SL_2(\s C)\sdr\s C^3$ of $\Gr$.

\section{Meromorphic functions}
\Proposition
Let $G$ be a connected complex linear-algebraic group with $G=G'$ and
an open subset $\Omega$ such that each element in $\Omega$ is semisimple.
Let $\Gamma$ be a Zariski-dense subgroup.

Then any $\Gamma$-invariant plurisubharmonic or meromorphic function is
constant and there exist no $\Gamma$-invariant hypersurface.
\endproclaim
\Proof
We may assume that $\Gamma$ is closed (in the Hausdorff topology).
Since $G=G'$, it follows that $H^0$ is a normal algebraic subgroup
for each Zariski-dense subgroup $H$. Therefore we may assume that
$\Gamma$ is discrete and furthermore it suffices (by induction on
$\dim(G)$ to demonstrate that the functions resp. hypersurfaces
are invariant under a positive-dimensional subgroup.
Now $G=G'$ implies that $\Gamma$ admits a finitely generated subgroup
$\Gamma_0$ which is still Zariski-dense. By the theorem of Selberg
$\Gamma_0$ admits a subgroup of finite index $\Gamma_1$ which is
torsion-free. Thus $\Gamma_1$ contains a semisimple element of infinite
order $\gamma$ which generates a subgroup $I$ whose Zariski-closure
$\bar I$ is a torus. $G=G'$ implies that this torus is contained in a
connected semisimple subgroup $S$ of $G$.
Now known results on subgroups in semisimple groups
\cite{HM}\cite{BKO} imply that the functions resp. hypersurfaces are
invariant under $\bar I$, which is positive-dimensional.
\qed

For this result it is essential to require $G=G'$ and not only
$G/G'$ reductive.
\Lemma
Let $G=\s C^*\times\s C^*$ and $\Gamma\simeq\s Z$
a (possibly Zariski-dense) discrete subgroup.

Then $G$ admits $\Gamma$-invariant
non-constant plurisubharmonic and meromorphic functions.
\endproclaim
\Proof
$G/\Gamma\simeq\s C^2/\Lambda$ with $\Lambda\simeq\s Z^3$.
Let $V=<\Lambda>_{\s R}$ the real subvector space of $\s C^2$
spanned by $\Lambda$
and $t:\s C^2\to\s C^2/V\simeq\s R$ a $\s R$-linear map.
Then $t^2$ yields a $\Gamma$-invariant plurisubharmonic
function on $G$.

Let $L=V\cap iV$ and $\gamma\in\Lambda\setminus L$. Let $H=<\gamma>_{\s C}$.
Then $H\ne L$, hence $H+L=\s C^2$. It follows that the $H$-orbits in
$G/\Gamma$ are closed and induce a fibration $G/\Gamma\to G/H\Gamma$ onto
a one-dimensional torus. One-dimensional tori are projective and therefore
admit non-constant meromorphic functions.
\qed

\References
\item{BO} Barth, W.; Otte, M.:
Invariante Holomorphe Funktionen auf reduktiven Liegruppen.
\sl Math. Ann. \bf 201\rm, 91--112 (1973)

\item{BKO} Berteloot, F.; Oeljeklaus, K.:
Invariant Plurisubharmonic Functions
and Hypersurfaces on Semisimple Complex Lie Groups.
\sl Math. Ann. \bf 281\rm, 513--530 (1988)

\item{HO} Huckleberry, A.T.; Oeljeklaus, E.:
On holomorphically separable complex solvmanifolds.
\sl Ann. Inst. Fourier \bf XXXVI 3\rm, 57--65 (1986)

\item{HM} Huckleberry, A.T.; Margulis, G.A.:
Invariant analytic hypersurfaces.
\sl Invent. Math. \bf 71\rm, 235--240 (1983)

\item{M1} Margulis, G.A.:
Free totally discontinuous groups of affine transformations.
\sl Soviet Math. Doklady \bf 28\rm, 435--439 (1983)

\item{M2} Margulis, G.A.:
Complete affine locally flat manifolds with free fundamental group.
\sl Zapiski Nau\v cn. Sem. Leningrad Otd. Mat. Inst. Steklov \bf 134\rm,
190-205 (1984) [in Russian]

\item{Mi} Milnor, J.:
On fundamental groups of completely affinely flat manifolds.
\sl Adv. Math. \bf 25\rm, 178--187 (1977)

\item{T} Tits, J.:
Free subgroups in linear groups.
\sl J. Algebra \bf 20\rm, 250--270 (1972)

\endarticle
%
% End of file
%
\\